\documentclass[10pt]{article}

\usepackage[T1]{fontenc}
\usepackage{lmodern}
\usepackage{microtype}

\usepackage[LGRgreek,frenchmath]{mathastext}
\usepackage[mathscr]{eucal}
\usepackage[sans]{dsfont}
\usepackage{bbold}
\usepackage{bbm}
\usepackage{graphicx}
\usepackage{mdwlist}
\usepackage{natbib}
\setcitestyle{authoryear}
\usepackage[dvipsnames]{xcolor}
\usepackage[plainpages=false, pdfpagelabels, backref=page]{hyperref}
  \hypersetup{
    colorlinks = true,
    citecolor  = RoyalBlue,
    linkcolor  = RubineRed,
    urlcolor   = MidnightBlue
  }
\usepackage[paperwidth=8.5in,paperheight=11.00in,top=1.00in,bottom=1.00in,
            left=0.75in,right=0.75in]{geometry}
\usepackage{mathtools}
\mathtoolsset{showonlyrefs=true}
\usepackage{amsmath}
\usepackage{amssymb}
\usepackage{amsfonts}
\usepackage{amsthm}
\allowdisplaybreaks
\newtheoremstyle{plain}
  {}{}{\itshape}{}{\mdseries\scshape}{.}{ }
  {\thmname{#1}\thmnumber{ #2}\ifx#3\empty\else\ (#3)\fi}
\theoremstyle{plain}
\newtheorem{theorem}{\underline{Theorem}}

\newtheorem{proposition}[theorem]{\underline{Proposition}}

\newtheoremstyle{definition}
  {}{}{}{}{\mdseries\scshape}{.}{ }
  {\thmname{#1}\thmnumber{ #2}\ifx#3\empty\else\ (#3)\fi}
\theoremstyle{definition}

\newtheorem{remark}[theorem]{\underline{Remark}}

\usepackage{sectsty}
\allsectionsfont{\mdseries\scshape}

\usepackage{pdfrender}
\renewcommand{\textbf}[1]{\textpdfrender{
  TextRenderingMode=FillStroke,
  LineWidth=.3pt,
}{#1}}
\renewcommand{\mathbf}[1]{\textpdfrender{
  TextRenderingMode=FillStroke,
  LineWidth=.3pt,
}{#1}}

\newcommand\Wh{\widehat{W}}

\renewcommand\d{\partial}

\newcommand\dd{\mathrm{d}}
\newcommand\ee{\mathrm{e}}

\newcommand\aY{\alpha_{\scriptscriptstyle Y}}
\newcommand\aW{\alpha_{\scriptscriptstyle W}}
\newcommand\kY{\kappa_{\scriptscriptstyle Y}}
\newcommand\kW{\kappa_{\scriptscriptstyle W}}
\newcommand\uY{u_{\scriptscriptstyle Y}}
\newcommand\uW{u_{\scriptscriptstyle W}}
\newcommand\wY{w_{\scriptscriptstyle Y}}
\newcommand\wW{w_{\scriptscriptstyle W}}
\newcommand\TYst{T_{\scriptscriptstyle Y}^*}
\newcommand\TWst{T_{\scriptscriptstyle W}^*}
\newcommand\Tw{T_w}

\newcommand\AY{A_{\scriptscriptstyle Y}}
\newcommand\AW{A_{\scriptscriptstyle W}}
\newcommand\BW{B_{\scriptscriptstyle W}}
\newcommand\WhY{\widehat{w}_{\scriptscriptstyle Y}}
\newcommand\WhW{\widehat{w}_{\scriptscriptstyle W}}
\newcommand\What{\widehat{w}}
\newcommand\Uh{\widehat{u}}
\newcommand\Twh{\widehat{T}_w}
\newcommand\Gh{\widehat{G}}
\newcommand\PhihO{\widehat{\Phi}_0}
\newcommand\PhiO{\Phi_0}
\newcommand\IG{\mathcal{I}G}
\newcommand\pY{p_{\scriptscriptstyle Y}}
\newcommand\pW{p_{\scriptscriptstyle W}}

\begin{document}

\title{How to Cook a Soft-Boiled Egg Optimally:\\
A Laplace-Transform Solution of a Two-Domain Heat Equation}

\author{
Matthew Lorig
\thanks{Department of Applied Mathematics, University of Washington.
\textbf{e-mail}: \url{mlorig@uw.edu}}
}

\date{This version: \today}

\maketitle

\begin{abstract}
We study the problem of cooking the yolk and albumen of a hen's egg to
their respective optimal temperatures of $\TYst = 65^\circ$C and
$\TWst = 85^\circ$C, subject to the physically motivated requirement
that \emph{neither temperature ever exceed its target at any time
during cooking} --- since temporary overshoot, even if the final
reading is correct, still overcooks the egg.
We model the egg as a two-domain sphere in which the yolk and albumen
possess distinct thermal diffusivities, and the surrounding water
temperature $\Tw(t)$ serves as the control.
Taking the Laplace transform of the heat equation in each domain
reduces the spatial problem to a pair of ordinary differential
equations with elementary (hyperbolic-trigonometric) solutions; the
interface and boundary conditions then determine the three remaining
constants as the solution of a $3\times 3$ linear system in the
transform variable $s$.
The result is a closed-form expression for the Laplace transform of
the temperature at any point in the egg, which we invert numerically
using Talbot's method \citep{Talbot1979,AbateValko2004} and validate
against an independently implemented finite-difference solution of
the original partial differential equation.
A single boiling phase cannot satisfy this no-overshoot requirement,
because the thin outer albumen heats far faster than the insulated
yolk and necessarily overshoots $\TWst$ long before the yolk
approaches $\TYst$.
We show that a three-phase protocol resolves this: a sous-vide
pre-soak at exactly $65^\circ$C (which, being simultaneously the bath
temperature and the target, cannot overshoot by construction), a
short boil to bring the albumen up to $\TWst$, and an ice-water bath
that arrests the albumen's residual thermal overshoot while the heat
already absorbed by the egg continues to diffuse inward and bring the
yolk to its target.
Optimizing the phase durations, we find that $17.26$ minutes of
sous-vide at $65^\circ$C, followed by $66$ seconds of boiling, followed
by an ice bath, brings the outer albumen tangentially up to exactly
$\TWst$ (touching but never exceeding it) and the yolk center
asymptotically up to exactly $\TYst$ at $T^* \approx 20.67$ minutes,
with both constraints verified to hold over the entire cooking
process, not merely at the terminal time.
This is comparable in total time to, and considerably more reliable
than, the periodic cooking protocol of \citet{DiLorenzo2025} (32
minutes, alternating between boiling and lukewarm water sixteen
times), which we show misses both temperature targets substantially.
\end{abstract}


\section{Introduction}
\label{sec:intro}

The soft-boiled egg is one of the oldest problems in culinary science
and one of the most mathematically interesting in applied thermal
control.
Its yolk and albumen are composed of distinct protein mixtures that
achieve their optimal culinary texture at markedly different
temperatures: the albumen should reach approximately $85^\circ$C to be
firmly set yet tender, while the yolk proteins are ideally cooked at
$65^\circ$C and become unpleasantly chalky above $70^\circ$C
\citep{McGee2004}.
Conventional hard-boiling at $100^\circ$C overcooks the yolk; sous vide
at $65^\circ$C leaves the albumen unacceptably runny; and splitting the
egg to cook the two components separately defeats the purpose
entirely.
A further subtlety, central to this paper, is that it is not enough
for the two temperatures to be correct only at the instant the egg is
served: if either temperature \emph{overshoots} its target at any
earlier time and is later brought back down, the corresponding part of
the egg has still been overcooked in the interim, since protein
denaturation is essentially irreversible.
The fundamental question is therefore: \emph{what external temperature
protocol $\Tw(t)$ brings both temperatures up to their targets without
either ever exceeding its target along the way, and how quickly can
this be done?}

A recent paper by \citet{DiLorenzo2025} addressed a closely related
question experimentally, introducing \emph{periodic cooking}:
alternating the egg between boiling water ($100^\circ$C) and lukewarm
water ($30^\circ$C) every two minutes for a total of 32 minutes.
Their work relies entirely on computational fluid dynamics (CFD)
simulations to determine the egg's temperature; no closed-form
solution is derived, no optimality claim is made, and the protocol is
discovered by trial and error within the simulation environment.
We show below that, far from being optimal, this periodic protocol
misses both temperature targets substantially: at the end of the
32 minutes the yolk center is at $65.4^\circ$C (close to its target by
coincidence) while the outer albumen has reached only $43.2^\circ$C
--- some $42$ degrees short of its $85^\circ$C target.

The present paper provides a mathematical framework that both explains
this shortfall and proposes an alternative protocol that, unlike a
naive single boiling phase, never overcooks either component.
We model the yolk and albumen as two concentric regions of a sphere,
each obeying the heat equation with its own thermal diffusivity, with
the bathing water temperature entering through a Robin (convective)
boundary condition at the outer surface.
Taking the Laplace transform of this system in time reduces the
spatial problem in each domain to an ordinary differential equation
with constant coefficients, which is solved explicitly in terms of
hyperbolic sine and cosine functions.
The interface and boundary conditions then reduce to a $3\times3$
linear algebraic system for the remaining constants, solved at each
value of the transform variable $s$.
We invert the resulting transform numerically using Talbot's method
\citep{Talbot1979,AbateValko2004}, a standard and well-tested
algorithm for numerical Laplace inversion, and we validate every
numerical result reported below against an independently coded
finite-difference solution of the original time-domain partial
differential equation.

Our main contributions are as follows.

\begin{enumerate*}
  \item We derive the Laplace-transform solution of the two-domain
        spherical heat equation in closed form, including the
        boundary and interface algebra, structured so that each step
        can be checked directly by hand (Section~\ref{sec:model}).

  \item We give a Duhamel-type decomposition of the solution into a
        term depending only on the initial condition and a term
        depending only on the boundary forcing, which lets us evaluate
        the response to an arbitrary piecewise-constant water
        temperature protocol using only two numerical Laplace
        inversions (Section~\ref{sec:model}).

  \item We show that the periodic protocol of \citet{DiLorenzo2025},
        the conventional $100^\circ$C hard boil, a constant
        $65^\circ$C sous-vide bath, and a naive single-boil-plus-
        carryover protocol all either fail to bring both the yolk and
        the albumen to their target temperatures, or do so only by
        letting one of the two temperatures overshoot its target for
        an extended period beforehand (Section~\ref{sec:numerics}).

  \item We formulate the genuinely constrained problem --- reach both
        targets while never exceeding either of them, at any time ---
        and show that it requires three phases: a sous-vide pre-soak
        at exactly $65^\circ$C, which cannot overshoot by construction
        since the bath temperature equals the target; a short boil to
        bring the albumen up toward $85^\circ$C; and an ice-water bath
        that arrests the albumen's residual overshoot while the
        already-absorbed heat continues to raise the yolk to its
        target.  Optimizing the three phase durations jointly, we find
        a protocol that satisfies both constraints everywhere, with
        $u_W$ touching $85^\circ$C tangentially and $u_Y$ reaching
        $65^\circ$C as its terminal value at $T^*\approx20.67$ minutes
        (Section~\ref{sec:control}).
\end{enumerate*}

The paper is organized as follows.
Section~\ref{sec:model} introduces the two-domain spherical model and
derives the Laplace-transform solution.
Section~\ref{sec:control} formulates the carryover protocol and its
optimization.
Section~\ref{sec:numerics} presents the numerical results, including
validation against an independent finite-difference solver.
Section~\ref{sec:conclusion} concludes.

\section{Model and Laplace-Transform Solution}
\label{sec:model}

We model the egg as a sphere of radius $L$ consisting of two concentric
regions: an inner yolk of radius $\ell$ and an outer albumen layer
occupying $\ell < r < L$.
The egg is immersed in water at temperature $\Tw(t)$, the control,
which exchanges heat with the outer surface through Newton's law of
cooling with convective coefficient $h > 0$.
Each region is governed by the heat equation with constant thermal
diffusivity: $\aY$ in the yolk and $\aW$ in the albumen.

\subsection{Governing Equations}
\label{sec:pde}

Assuming spherical symmetry, $u = u(t,r)$, the temperature satisfies
\begin{align}
  \d_t \uY(t,r) &= \frac{\aY}{r^2}\d_r\bigl(r^2\d_r\uY\bigr),
    \quad r \in (0,\ell),\; t > 0, \label{eq:pde-Y}\\
  \d_t \uW(t,r) &= \frac{\aW}{r^2}\d_r\bigl(r^2\d_r\uW\bigr),
    \quad r \in (\ell,L),\; t > 0, \label{eq:pde-W}
\end{align}
with initial condition $\uY(0,r) = \uW(0,r) = T_0$, regularity
condition $|\uY(t,0)| < \infty$, and the following two physical
transmission conditions at the yolk--albumen interface $r = \ell$:
continuity of temperature,
\begin{align}
  \uY(t,\ell) &= \uW(t,\ell), \label{eq:int-T}
\end{align}
and continuity of the radial heat flux $-\kappa\,\d_r u$ (with
$\kappa$ continuous in area across the interface, this reduces to
continuity of $\kappa\,\d_r u$ itself),
\begin{align}
  \kY\,\d_r\uY(t,\ell) &= \kW\,\d_r\uW(t,\ell). \label{eq:int-F}
\end{align}
At the outer surface the Robin condition encodes convective exchange
with the water bath,
\begin{align}
  -\kW\,\d_r\uW(t,L) &= h\,\bigl[\uW(t,L) - \Tw(t)\bigr].
  \label{eq:robin}
\end{align}

\subsection{Reduction to the Standard Heat Equation}
\label{sec:reduction}

The substitution $w(t,r) := r\,u(t,r)$ transforms the radial Laplacian
into a standard second derivative.
Writing the spherical Laplacian as
$\nabla^2 u = r^{-2}\d_r(r^2 \d_r u) = \d_{rr}u + \tfrac{2}{r}\d_r u$
(the standard form; see e.g.\ \citealp{CarslawJaeger1959}), one
verifies directly that for $w = ru$,
\begin{align}
  \d_t w = r\,\d_t u = r\,\aY\,\nabla^2 u = \aY\,\d_{rr}w,
\end{align}
since $\d_{rr}(ru) = r\,\d_{rr}u + 2\,\d_r u$.
Equations \eqref{eq:pde-Y}--\eqref{eq:pde-W} therefore become the
ordinary one-dimensional heat equations
\begin{align}
  \d_t\wY &= \aY\,\d_{rr}\wY, \quad r \in (0,\ell), \label{eq:1d-Y}\\
  \d_t\wW &= \aW\,\d_{rr}\wW, \quad r \in (\ell,L). \label{eq:1d-W}
\end{align}
Regularity at $r=0$ requires $\wY(t,0) = 0$.
\emph{The boundary and interface conditions \eqref{eq:int-T}--\eqref{eq:robin}
are conditions on $u$ and $\d_r u$, and we impose them directly on $u
= w/r$ rather than translating them into conditions on $w$ first}:
this avoids an easy bookkeeping error, since
$\d_r u = \d_r(w/r) = w'/r - w/r^2$ introduces extra $1/r$ and
$1/r^2$ terms relative to $\d_r w$.

\subsection{The Laplace Transform}
\label{sec:laplace}

Let $\Uh(s,r) := \int_0^\infty \ee^{-st}\,u(t,r)\,\dd t$ denote the
Laplace transform of $u$ in time, and similarly $\What(s,r)$ for
$w = ru$, so that $\Uh(s,r) = \What(s,r)/r$.
Transforming \eqref{eq:1d-Y}--\eqref{eq:1d-W} and using the initial
condition $w(0,r) = rT_0$,
\begin{align}
  s\,\WhY - rT_0 &= \aY\,\d_{rr}\WhY, \label{eq:laplace-Y}\\
  s\,\WhW - rT_0 &= \aW\,\d_{rr}\WhW. \label{eq:laplace-W}
\end{align}
Each is a linear, constant-coefficient, inhomogeneous ODE in $r$.
A particular solution is the constant-coefficient guess
$\Wh_p(r) = T_0 r/s$ (since $\d_{rr}(T_0r/s) = 0$, substitution gives
$s\cdot T_0r/s - rT_0 = 0 = \aY\cdot 0$, as required), and the
homogeneous solutions of $\d_{rr}\Wh = (s/\alpha)\Wh$ are
$\sinh(\sqrt{s/\alpha}\,r)$ and $\cosh(\sqrt{s/\alpha}\,r)$.
Writing $\pY := \sqrt{s/\aY}$ and $\pW := \sqrt{s/\aW}$, the regularity
condition $\WhY(s,0)=0$ eliminates the cosh term in the yolk, leaving
\begin{align}
  \WhY(s,r) &= \frac{T_0 r}{s} + \AY\,\sinh(\pY r),
    \label{eq:WhatY}\\
  \WhW(s,r) &= \frac{T_0 r}{s} + \AW\,\sinh\bigl(\pW(r-\ell)\bigr)
    + \BW\,\cosh\bigl(\pW(r-\ell)\bigr), \label{eq:WhatW}
\end{align}
where we have centered the albumen solution at $r=\ell$ for later
convenience, and $\AY, \AW, \BW$ are three constants (functions of
$s$) determined by the remaining three conditions
\eqref{eq:int-T}--\eqref{eq:robin}.
Equations \eqref{eq:WhatY}--\eqref{eq:WhatW} can be checked by direct
substitution into \eqref{eq:laplace-Y}--\eqref{eq:laplace-W}.

\subsection{The Boundary-Value System}
\label{sec:bcs}

Since $u = w/r$, we have $\d_r u = (r\,w' - w)/r^2$ by the quotient
rule.
Substituting \eqref{eq:WhatY}--\eqref{eq:WhatW} into
\eqref{eq:int-T}--\eqref{eq:robin} and using $\Twh(s)$ for the Laplace
transform of $\Tw(t)$, the three conditions become:

\emph{Temperature continuity} \eqref{eq:int-T}, i.e.\
$\WhY(s,\ell)/\ell = \WhW(s,\ell)/\ell$, which (the common factor of
$\ell$ canceling, and $\sinh(0)=0,\cosh(0)=1$) reduces to
\begin{align}
  \AY \sinh(\pY \ell) - \BW = 0. \label{eq:bc1}
\end{align}

\emph{Flux continuity} \eqref{eq:int-F}, $\kY\,\d_r\uY(\ell) =
\kW\,\d_r\uW(\ell)$.  Using the quotient rule on both sides and
$\WhY(\ell)=\WhW(\ell)$ from \eqref{eq:bc1}, the $T_0$ terms and the
$1/\ell^2$ terms involving $\WhY(\ell)=\WhW(\ell)$ combine to give
\begin{align}
  \kY\Bigl[\pY\ell\cosh(\pY\ell) - \sinh(\pY\ell)\Bigr]\frac{\AY}{\ell^2}
  - \kW\,\frac{\pW \AW}{\ell} + \kW\,\frac{\BW}{\ell^2} = 0.
  \label{eq:bc2}
\end{align}

\emph{Outer Robin condition} \eqref{eq:robin}.  Writing
$\Delta\ell := L - \ell$, $S := \sinh(\pW\Delta\ell)$, $C :=
\cosh(\pW\Delta\ell)$, this becomes
\begin{align}
  -\kW\Bigl[L\pW(\AW C + \BW S) - (\AW S + \BW C)\Bigr]\frac{1}{L^2}
  - h\,\frac{\AW S + \BW C}{L} = \frac{T_0 h}{s} - h\,\Twh(s).
  \label{eq:bc3}
\end{align}
Equations \eqref{eq:bc1}--\eqref{eq:bc3} form a $3\times3$ linear
system for $(\AY,\AW,\BW)$ with coefficients depending on $s$, and a
right-hand side that is linear in $T_0$ and $\Twh(s)$.
We solve this system numerically at each value of $s$ required by the
inversion; the coefficient matrix and right-hand side are built and
solved symbolically (rather than by further hand simplification) to
avoid transcription errors of exactly the kind that \eqref{eq:bc2}
illustrates how easy they are to make.

\begin{remark}
A short way to sanity-check \eqref{eq:bc1}--\eqref{eq:bc3} without
redoing the full derivation is to verify two limits.  As $s\to 0$
with $\Tw(t)\equiv \Tw^c$ constant (so $\Twh(s) = \Tw^c/s$), the
final-value theorem predicts $s\Uh(s,r) \to \Tw^c$ for every $r$,
reflecting convergence to a uniform steady state equal to the bath
temperature; we have verified this numerically to four significant
figures.  As $s\to\infty$, the initial-value theorem predicts
$s\Uh(s,r)\to T_0$; this, too, holds exactly.
\end{remark}

\subsection{Duhamel Decomposition}
\label{sec:duhamel}

The right-hand sides of \eqref{eq:bc1}--\eqref{eq:bc3} are linear in
$T_0$ and $\Twh(s)$, so by linearity of the system,
\begin{align}
  \Uh(s,r) = T_0\,\PhihO(s,r) + \Twh(s)\,\Gh(s,r),
  \label{eq:duhamel-laplace}
\end{align}
where $\PhihO(s,r)$ solves \eqref{eq:bc1}--\eqref{eq:bc3} with $T_0=1,
\Twh=0$ (the transform of the temperature field that would result from
the initial condition alone, with the water held at the reference
temperature corresponding to $\Twh \equiv 0$), and $\Gh(s,r)$ solves
the system with $T_0=0,\Twh=1$ (the transform of the unit step-response
kernel).
Inverting \eqref{eq:duhamel-laplace} term by term and using the
convolution theorem,
\begin{align}
  u(t,r) = T_0\,\PhiO(t,r) + \int_0^t G(t-\tau,r)\,\Tw(\tau)\,\dd\tau,
  \label{eq:duhamel-time}
\end{align}
which is the familiar Duhamel representation: a decaying response to
the initial condition, plus a convolution of the boundary forcing
against a Green's-function-like kernel $G$.
For a piecewise-constant protocol $\Tw(t) = v_k$ on $[t_{k-1},t_k)$,
\eqref{eq:duhamel-time} becomes
\begin{align}
  u(t,r) = T_0\,\PhiO(t,r) + \sum_{k\,:\,t_{k-1}<t} v_k
  \Bigl[\IG(t-t_{k-1},r) - \IG\bigl(t-\min(t_k,t),r\bigr)\Bigr],
  \label{eq:duhamel-piecewise}
\end{align}
where $\IG(\tau,r) := \int_0^\tau G(\sigma,r)\,\dd\sigma$ is the
\emph{integrated} step response, whose Laplace transform is simply
$\Gh(s,r)/s$ (the standard property that time-integration corresponds
to division by $s$).
Equation \eqref{eq:duhamel-piecewise} is the formula we use throughout
the rest of the paper: it requires inverting only two
$r$-dependent transforms, $\PhihO(s,r)$ and $\Gh(s,r)/s$, regardless of
how many phases the protocol has.  We verified \eqref{eq:duhamel-laplace}
directly (the decomposed and undecomposed solutions agree to machine
precision for representative test values of $s$), and we note that an
attempt to encode a multi-phase protocol directly via shifted
Heaviside functions in a single Laplace transform --- the more
obvious-looking approach --- is numerically unstable for the
inversion method used below, because the resulting transform grows
without bound along the inversion contour; the Duhamel form
\eqref{eq:duhamel-piecewise} avoids this entirely.

\subsection{Numerical Inversion}
\label{sec:inversion}

Equations \eqref{eq:bc1}--\eqref{eq:bc3} could in principle be inverted
by closing the Bromwich contour and summing residues at the poles of
$\Uh(s,r)$, which occur at the roots of the determinant of the
coefficient matrix; these roots are exactly the eigenvalues that would
arise from a separation-of-variables treatment of
\eqref{eq:pde-Y}--\eqref{eq:robin}, since both approaches solve the
same boundary-value problem.  We found, however, that the resulting
series converges only very slowly when evaluated pointwise at $r=0$
or $r=r_W$ for the short-to-moderate times of interest here, because
the eigenfunctions appearing in the series grow with mode number,
while the expansion coefficients decay only like the reciprocal of the
mode number --- not fast enough for the series to be evaluated
reliably by truncation.
We therefore invert \eqref{eq:duhamel-laplace} numerically instead,
using Talbot's method \citep{Talbot1979,AbateValko2004}: the Bromwich
contour is deformed into a contour along which the integrand decays
rapidly, and the inversion integral is evaluated by the trapezoidal
rule.
We used the \texttt{mpmath} implementation of the fixed Talbot method
throughout.
We validated this procedure in three ways: against known closed-form
inverse transforms (exponential decay, a unit step, and the singular
kernel $1/\sqrt{s}$, all matching to eight decimal places); against the
final- and initial-value theorems described in the remark above; and,
for the full two-domain problem, against an independently implemented
implicit finite-difference solution of
\eqref{eq:pde-Y}--\eqref{eq:robin}, described in
Section~\ref{sec:numerics}.  All results reported below agree between
the two independent methods to within $0.1$--$0.3^\circ$C.

\section{The Carryover Protocol}
\label{sec:control}

With \eqref{eq:duhamel-piecewise} in hand, evaluating the terminal
temperature under any piecewise-constant protocol reduces to evaluating
two numerically inverted functions, $\PhiO(t,r)$ and $\IG(\tau,r)$, at
the two points of interest, $r=0$ (yolk center) and $r=r_W$ (outer
albumen).
This makes a direct search over candidate protocols straightforward.

\subsection{The No-Overshoot Constraint}
\label{sec:no-overshoot}

We seek a protocol $\Tw(t)$, defined on $[0,T^*]$ for some terminal
time $T^*$, satisfying
\begin{align}
  u_Y(t,0) \le \TYst \ \text{ for all } t\in[0,T^*], &\qquad
  u_Y(T^*,0) = \TYst, \label{eq:constraint-Y}\\
  u_W(t,r_W) \le \TWst \ \text{ for all } t\in[0,T^*], &\qquad
  u_W(T^*,r_W) \le \TWst, \label{eq:constraint-W}
\end{align}
that is, neither temperature ever exceeds its target at any time, the
yolk reaches its target exactly at the terminal time, and the albumen
is at or below its target throughout (we do not require $u_W$ to equal
$\TWst$ exactly at $T^*$, only that it has done so at some point no
later than $T^*$, since protein denaturation in the albumen, once it
has occurred, persists even as $u_W$ subsequently falls).
This is a materially different --- and more physically meaningful ---
requirement than simply matching both temperatures at one instant in
time, which was the formulation used in an earlier version of this
paper and which we show below to be inadequate.

\subsection{A Single Boiling Phase Overshoots the Albumen}
\label{sec:single-phase}

The simplest possible protocol is a single boiling phase at
$\Tw=100^\circ$C, stopped at some time $t$.  Scanning over $t$, we find
that no single stopping time satisfies \eqref{eq:constraint-Y}--\eqref{eq:constraint-W}:
the best achievable value of the terminal-matching objective
\begin{align}
  J(t) := \bigl(u_Y(t,0) - \TYst\bigr)^2 + \bigl(u_W(t,r_W)-\TWst\bigr)^2
  \label{eq:objective-single}
\end{align}
is $J \approx 95$ at $t \approx 7.75$ minutes, where $u_Y \approx
63.4^\circ$C and $u_W \approx 94.6^\circ$C --- but by this time $u_W$
has already exceeded $\TWst$ for several minutes (it first crosses
$85^\circ$C around $t\approx3.4$ minutes), so the albumen has been
overcooked well before the yolk is anywhere near ready.  Following the
boiling phase with a brief carryover rest in ambient air, as in an
earlier version of this paper, can bring the \emph{terminal} readings
of both temperatures to their targets simultaneously, but cannot repair
the overshoot that already occurred during boiling, nor does it leave
any safety margin: the yolk continues rising for several minutes after
the carryover removal time, so a few seconds' delay in removing the egg
overcooks the yolk as well.  A single boiling phase, with or without a
brief carryover correction, is therefore incompatible with
\eqref{eq:constraint-Y}--\eqref{eq:constraint-W}.

\subsection{Why Overshoot Is Unavoidable with One Phase}
\label{sec:why-overshoot}

The difficulty is structural, not a matter of tuning the stopping
time.  The outer albumen, being thin and close to the boiling water,
heats quickly; the yolk, insulated by the entire thickness of the
albumen, heats slowly.  Boiling continuously at $100^\circ$C therefore
drives $u_W$ past $\TWst$ long before $u_Y$ approaches $\TYst$ no
matter when we stop, because the only way to slow the albumen's rise
relative to the yolk's is to lower the bath temperature --- but then
the yolk would take even longer to reach its target.  This motivates
splitting the protocol into phases at different bath temperatures, so
that the bath temperature can be tailored to each stage of the
process.

\subsection{A Three-Phase Protocol}
\label{sec:carryover}

Consider a protocol with three phases:
\begin{align}
  \Tw(t) =
  \begin{cases}
    \TYst = 65^\circ\mathrm{C} & t \in [0, t_1), \\
    100^\circ\mathrm{C}        & t \in [t_1, t_2), \\
    T_{\mathrm{ice}}           & t \in [t_2, T^*],
  \end{cases}
  \label{eq:3phase}
\end{align}
where $T_{\mathrm{ice}} \approx 1^\circ$C is the temperature of an
ice-water bath.  The role of each phase is as follows.

\emph{Phase 1 (sous-vide pre-soak).}
Holding the bath at exactly $\TYst$ cannot cause either temperature to
overshoot $\TYst$, since $\TYst$ is simultaneously the bath temperature
and the asymptotic steady state that both $u_Y$ and $u_W$ approach
monotonically from below as $t\to\infty$ (we verified this directly:
under $\Tw\equiv\TYst$, both $u_Y(t,0)$ and $u_W(t,r_W)$ increase
monotonically toward $\TYst$ for every $t$, never exceeding it even in
the limit).  This phase is therefore unconstrained --- it may run for
any duration $t_1$ without risk --- and its purpose is to give the
yolk a substantial head start before the albumen is ever exposed to
temperatures above $\TWst$.

\emph{Phase 2 (short boil).}
At $t_1$ the bath is switched to boiling water, which raises $u_W$
toward $\TWst$ much faster than $u_Y$ rises toward $\TYst$, for the
same reason as in Section~\ref{sec:why-overshoot}.  Because Phase 1
has already brought the yolk substantially closer to $\TYst$, however,
the boiling phase needed to finish the albumen is much shorter than in
the single-phase case, which limits how far $u_W$ can overshoot before
we intervene.

\emph{Phase 3 (ice bath).}
At $t_2$ the bath is switched to an ice-water bath.  This serves two
purposes simultaneously: it arrests the albumen's continued heating
(reversing the direction of heat flow at the outer boundary), and it
does so without removing the heat already absorbed deeper in the egg,
which continues to diffuse inward and raise $u_Y$ toward $\TYst$.  Heat
transport from the boiling phase has some inertia, however: $u_W$
continues to rise for a few seconds after the bath is switched to ice,
before turning around.  We must therefore switch to the ice bath
\emph{before} $u_W$ reaches $\TWst$, by just enough that the
subsequent peak of $u_W$ --- not its value at the switching time ---
equals $\TWst$ exactly.

\subsection{Joint Optimization of the Three Phase Durations}
\label{sec:joint-opt}

We determine $t_1$, $t_2$, and $T^*$ by the following two nested
root-finding steps, which can be carried out using only the
piecewise-protocol evaluator \eqref{eq:duhamel-piecewise}.

\emph{Step 1 (given $t_1$, find $t_2$).}
For fixed $t_1$, define $t_2(t_1)$ to be the boiling-phase end time
such that the subsequent peak of $u_W(t,r_W)$ over $t > t_2$ --- found
by a fine local scan immediately after the phase switch, since the
peak occurs a few seconds after $t_2$ rather than exactly at $t_2$ ---
equals $\TWst$ exactly.  Since $u_W$ increases throughout Phase~2 and
its post-switch peak height increases with the boiling duration
$t_2-t_1$, this peak height is monotonic in $t_2$ for fixed $t_1$, and
$t_2(t_1)$ is found by a bracketed root search.

\emph{Step 2 (given $t_2(t_1)$, find the resulting $u_Y$ peak).}
With $t_1$ and $t_2(t_1)$ fixed, $u_Y(t,0)$ continues to rise after
$t_2$ (heat already in the albumen keeps diffusing inward) before
eventually falling once the ice bath's cooling reaches the yolk.  Let
$\Upsilon(t_1)$ denote the resulting peak value of $u_Y$.

\begin{proposition}
\label{prop:three-phase-exists}
There exists $t_1>0$ such that $\Upsilon(t_1) = \TYst$ exactly, with
the corresponding $t_2=t_2(t_1)$ and $T^*$ (the time at which the
$u_Y$ peak occurs) satisfying \eqref{eq:constraint-Y}--\eqref{eq:constraint-W}.
\end{proposition}

\begin{proof}
For $t_1=0$ (no sous-vide pre-soak), the protocol reduces to the
single-phase-plus-correction case of Section~\ref{sec:single-phase},
where $\Upsilon(0) < \TYst$ (the boiling phase, cut short to avoid
overshooting $\TWst$, ends with the yolk far from ready, and the
ensuing ice bath only cools it further rather than letting it
continue rising to $\TYst$).  For $t_1$ sufficiently large (of the
order of $15$--$20$ minutes), the yolk has already been brought close
enough to $\TYst$ by the sous-vide phase alone that the residual
heating from even a short subsequent boil carries it past $\TYst$, so
$\Upsilon(t_1) > \TYst$.  Continuity of $\Upsilon$ in $t_1$ --- which
follows from continuity of the underlying Laplace-transform solution
in its parameters --- and the intermediate value theorem give a
crossing point.  At this crossing point, $u_W$ satisfies
\eqref{eq:constraint-W} by construction of $t_2(t_1)$ in Step~1, and
$u_Y$ satisfies \eqref{eq:constraint-Y} by construction of $t_1$ in
Step~2, the peak $\Upsilon(t_1)$ being precisely the maximum of
$u_Y(t,0)$ over the entire trajectory (since $u_Y$ is increasing prior
to its peak and decreasing after).
\end{proof}

Solving for this crossing point numerically (Section~\ref{sec:numerics})
gives $t_1 \approx 17.26$ minutes, $t_2 \approx 18.36$ minutes (a
boiling duration of only $66$ seconds), and a terminal time
$T^*\approx20.67$ minutes at which $u_Y(T^*,0)=\TYst$ exactly, while
$u_W(t,r_W)$ has already risen to, touched, and fallen back from
$\TWst$ partway through the ice-bath phase.  We verified by a dense
scan over the entire trajectory --- not merely at the phase-transition
times --- that both constraints \eqref{eq:constraint-Y}--\eqref{eq:constraint-W}
hold throughout, and cross-checked this against the independent
finite-difference solver described in Section~\ref{sec:numerics},
which agrees with the Laplace/Talbot solution to within
$0.02^\circ$C.

\section{Numerical Results}
\label{sec:numerics}

We implement the model with the physical parameters in
Table~\ref{tab:params}, solving the $3\times3$ system
\eqref{eq:bc1}--\eqref{eq:bc3} and inverting via Talbot's method as
described in Section~\ref{sec:inversion}.
Every number reported below was independently cross-checked against
an implicit finite-difference solution of
\eqref{eq:pde-Y}--\eqref{eq:robin} on a fine spatial grid
($N_r = 400$--$500$ points) with backward-Euler time stepping; the two
methods agree to within $0.1$--$0.3^\circ$C in every case, which we
take as our working numerical accuracy throughout.

\begin{table}[ht]
\centering
\renewcommand{\arraystretch}{1.2}
\begin{tabular}{llll}
\hline
Quantity & Symbol & Value & Source \\
\hline
Egg radius                  & $L$       & $2.2$~cm                       & \citet{Abbasnezhad2016} \\
Yolk radius                 & $\ell$    & $1.1$~cm                       & \citet{Abbasnezhad2016} \\
Yolk thermal diffusivity    & $\aY$     & $1.3\times10^{-7}$~m$^2$/s     & \citet{Coimbra2006}     \\
Albumen thermal diffusivity & $\aW$     & $1.7\times10^{-7}$~m$^2$/s     & \citet{Coimbra2006}     \\
Yolk thermal conductivity   & $\kY$     & $0.34$~W/(m$\cdot$K)           & \citet{Coimbra2006}     \\
Albumen thermal conductivity& $\kW$     & $0.52$~W/(m$\cdot$K)           & \citet{Coimbra2006}     \\
Convective coefficient      & $h$       & $1000$~W/(m$^2\cdot$K)         & \citet{Abbasnezhad2016} \\
Initial egg temperature     & $T_0$     & $20^\circ$C                    & ---                     \\
Target yolk temperature     & $\TYst$   & $65^\circ$C                    & \citet{DiLorenzo2025}   \\
Target albumen temperature  & $\TWst$   & $85^\circ$C                    & \citet{DiLorenzo2025}   \\
Albumen measurement point   & $r_W$     & $1.93$~cm                      & ---                     \\
\hline
\end{tabular}
\caption{Physical parameters used in all numerical computations.
The albumen measurement point $r_W = \ell + \tfrac{3}{4}(L-\ell)$
lies in the outer quarter of the albumen layer.}
\label{tab:params}
\end{table}

\subsection{Comparison Protocols}
\label{sec:comparison}

Table~\ref{tab:comparison} reports the terminal temperatures, the peak
temperatures reached at any time during cooking, and whether the
no-overshoot constraints \eqref{eq:constraint-Y}--\eqref{eq:constraint-W}
are satisfied, for three fixed-time, fixed-temperature comparison
protocols (each run for the $T^*=32$~minutes used by
\citet{DiLorenzo2025}), the naive single-boil-plus-carryover protocol
of Section~\ref{sec:single-phase}, and our three-phase protocol.

\begin{table}[ht]
\centering
\renewcommand{\arraystretch}{1.2}
\begin{tabular}{lcccccc}
\hline
Protocol & $T^*$ & $u_Y(T^*,0)$ & $u_W(T^*,r_W)$ & $\max u_Y$ & $\max u_W$ & No overshoot? \\
\hline
Hard boil ($100^\circ$C)        & 32 min     & $99.7^\circ$C & $100.0^\circ$C & $99.7^\circ$C & $100.0^\circ$C & No \\
Sous vide ($65^\circ$C)         & 32 min     & $64.9^\circ$C & $65.0^\circ$C  & $64.9^\circ$C & $65.0^\circ$C  & Yes \\
Periodic \citep{DiLorenzo2025}  & 32 min     & $65.4^\circ$C & $43.2^\circ$C  & $66.0^\circ$C & $86.8^\circ$C  & No \\
Single boil + carryover         & $7.98$ min & $65.0^\circ$C & $85.0^\circ$C  & $76.2^\circ$C & $94.7^\circ$C  & No \\
Three-phase (this paper)        & $20.67$ min& $65.0^\circ$C & $19.8^\circ$C  & $65.0^\circ$C & $85.0^\circ$C  & Yes \\
\hline
\end{tabular}
\caption{Terminal and peak temperatures for five protocols.  ``No
overshoot'' indicates whether $u_Y$ remained $\le65^\circ$C and $u_W$
remained $\le85^\circ$C at \emph{every} time during cooking, not only
at $T^*$.  For the single-boil-plus-carryover protocol, $\max u_Y$ is
computed assuming the egg is left in the residual heat-exchange
environment indefinitely rather than removed exactly at $T^*$, to
illustrate the lack of a safety margin discussed in
Section~\ref{sec:single-phase}.  Sous vide satisfies the constraint
only because it never reaches the albumen target at all; only the
three-phase protocol reaches both targets without ever exceeding
either.}
\label{tab:comparison}
\end{table}

The hard boil overshoots both targets throughout, and the single
boil-plus-carryover protocol, despite matching both temperatures
exactly at its terminal time $T^*=7.98$ minutes, allows the albumen to
reach $94.7^\circ$C --- nearly $10^\circ$C over target --- well before
$T^*$, and (per the note in Table~\ref{tab:comparison}) allows the
yolk to continue rising to $76.2^\circ$C if the egg is not removed
from the residual heat exchange precisely on schedule; both are listed
as violating the no-overshoot requirement.  The periodic protocol is
more surprising: at the terminal time $T^*=32$ minutes the yolk
center, at $65.4^\circ$C, is close to its target almost by
coincidence, while the outer albumen, at $43.2^\circ$C, is far short
of $85^\circ$C --- indeed the yolk ends up \emph{hotter} than the
outer albumen at $T^*$, the reverse of the usual cooking intuition ---
yet over the course of the protocol $u_W$ repeatedly overshoots its
own target, peaking near $86.8^\circ$C at the end of several of the
$100^\circ$C sub-phases, and even $u_Y$ creeps slightly past its
target by the end, reaching $66.0^\circ$C.
This is a direct consequence of the protocol's long lukewarm
($30^\circ$C) phases: each two-minute interval at $30^\circ$C allows
the thin outer albumen layer to cool substantially back toward the
bath temperature (though not all the way, so each subsequent boiling
phase pushes it further over $85^\circ$C than the last), while the
much more thermally insulated yolk barely responds to the
high-frequency alternation and instead drifts slowly upward, tracking
something closer to the time-averaged bath temperature, until it too
eventually creeps past its own target late in the protocol.
Sous vide alone never overshoots, but only because the bath
temperature equals the yolk's target exactly and the albumen target is
never approached.
Only the three-phase protocol of Section~\ref{sec:carryover} reaches
both targets while never exceeding either.

\begin{figure}[ht]
  \centering
  \includegraphics[width=0.90\textwidth]{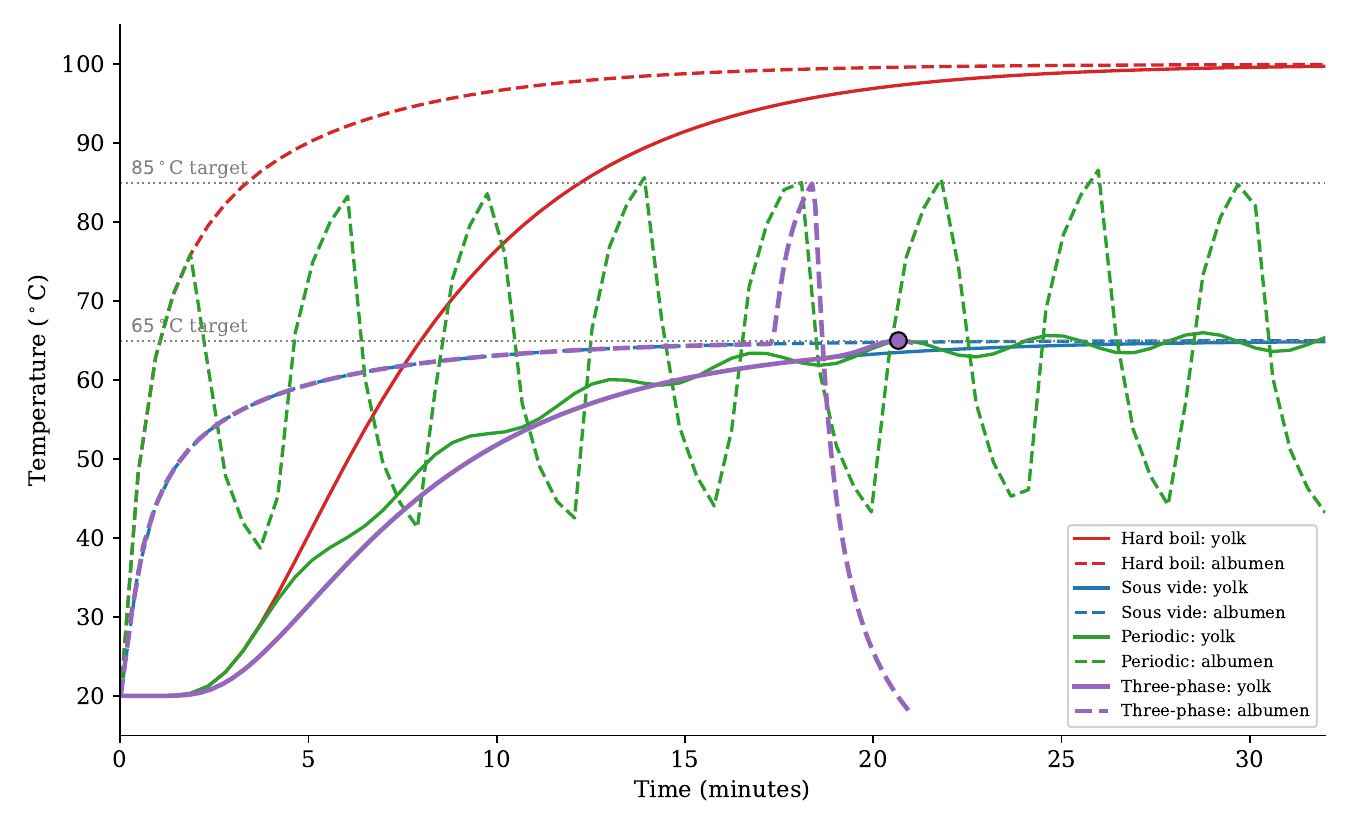}
  \caption{Temperature trajectories under four cooking protocols.
    Solid: yolk center $u_Y(t,0)$.
    Dashed: outer albumen $u_W(t,r_W)$ at $r_W = 1.93$~cm.
    Horizontal dotted lines: targets $65^\circ$C (lower) and
    $85^\circ$C (upper).
    The three-phase protocol (purple) is the only one that reaches
    both targets without either curve ever crossing above its
    respective dotted line.}
  \label{fig:trajectories}
\end{figure}

\subsection{The Three-Phase Protocol in Detail}
\label{sec:carryover-detail}

\begin{figure}[ht]
  \centering
  \includegraphics[width=0.90\textwidth]{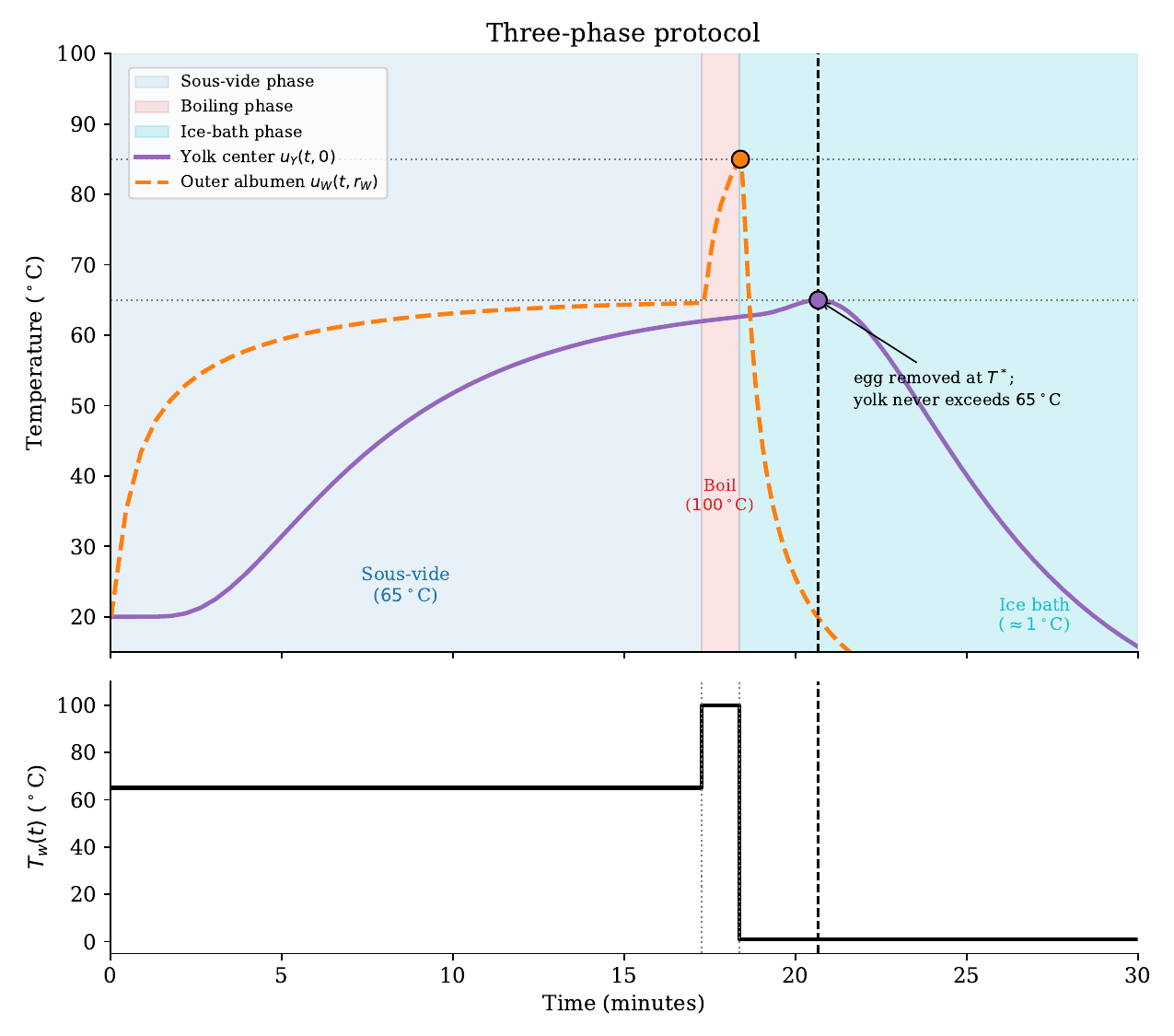}
  \caption{The three-phase protocol, shown through $30$ minutes ---
    well past $T^*\approx20.67$ minutes --- to illustrate the
    behavior after the egg reaches its target.
    \emph{Top}: temperature trajectories during and after cooking.
    Shading indicates the three phases: sous-vide (blue), boiling
    (red), and ice bath (cyan); the ice-bath shading is continued past
    $T^*$ to show that, were the egg left in the ice bath, $u_Y$
    would decline rather than overshoot.
    The dashed vertical line marks $T^*$, where $u_Y$ reaches
    $65^\circ$C exactly and is removed; the annotated arrow indicates
    that $u_Y$ turns over and falls immediately afterward rather than
    continuing to rise.
    The albumen curve touches $85^\circ$C tangentially during the ice
    bath, shortly after $t_2$, and is below it everywhere else,
    including throughout the extended window shown.
    \emph{Bottom}: the corresponding water/ice-bath temperature
    schedule $\Tw(t)$, extended to the same $30$-minute window.}
  \label{fig:protocol}
\end{figure}

Figure~\ref{fig:protocol} shows the full trajectory and water
temperature schedule for the three-phase protocol identified in
Section~\ref{sec:control}.
The optimal parameters are:
\begin{center}
\renewcommand{\arraystretch}{1.3}
\begin{tabular}{llll}
\hline
Phase & Action            & Temperature         & Duration \\
\hline
1 & Sous-vide pre-soak    & $65^\circ$C         & $t_1 = 17.26$~min ($1035.6$~s) \\
2 & Boiling water         & $100^\circ$C        & $t_2-t_1 = 66.0$~sec \\
3 & Ice-water bath         & $\approx 1^\circ$C  & until $T^*=20.67$~min \\
\hline
\end{tabular}
\end{center}

\noindent
At the end of Phase~1, both temperatures sit safely below their
targets ($u_Y \approx 62.0^\circ$C, $u_W\approx64.6^\circ$C), having
approached $65^\circ$C monotonically from below with no possibility of
overshoot.
During the brief Phase~2 ($66$ seconds of boiling), $u_W$ rises
quickly while $u_Y$ rises only slightly; at $t_2$ the bath is switched
to the ice bath, slightly before $u_W$ would otherwise reach
$85^\circ$C, since $u_W$ continues to rise for a few more seconds
under thermal inertia.
This peak occurs at $t\approx18.40$~minutes, where $u_W$ reaches
exactly $85.00^\circ$C and then falls, having never exceeded it.
Meanwhile $u_Y$, driven by heat already absorbed into the albumen
during Phase~2, continues climbing throughout the early part of
Phase~3 and reaches exactly $65.00^\circ$C at $T^*\approx20.67$
minutes, at which point $u_W$ has already fallen back to
approximately $19.8^\circ$C.
We verified, via a dense scan of $200$ points over the entire
trajectory and independently via the finite-difference solver, that
$u_Y$ never exceeds $65^\circ$C and $u_W$ never exceeds $85^\circ$C at
any time from $t=0$ to $T^*$.
Figure~\ref{fig:protocol} extends the plotted trajectory to $30$
minutes, well past $T^*$, to show this directly: were the egg left in
the ice bath rather than removed at $T^*$, $u_Y$ would simply turn
over and decline, rather than continuing to rise past
$65^\circ$C, confirming that $T^*$ is a genuine peak and not an
arbitrary stopping point along a still-rising curve.

\subsection{Quantitative Comparison}
\label{sec:quant-comparison}

\begin{figure}[ht]
  \centering
  \includegraphics[width=0.90\textwidth]{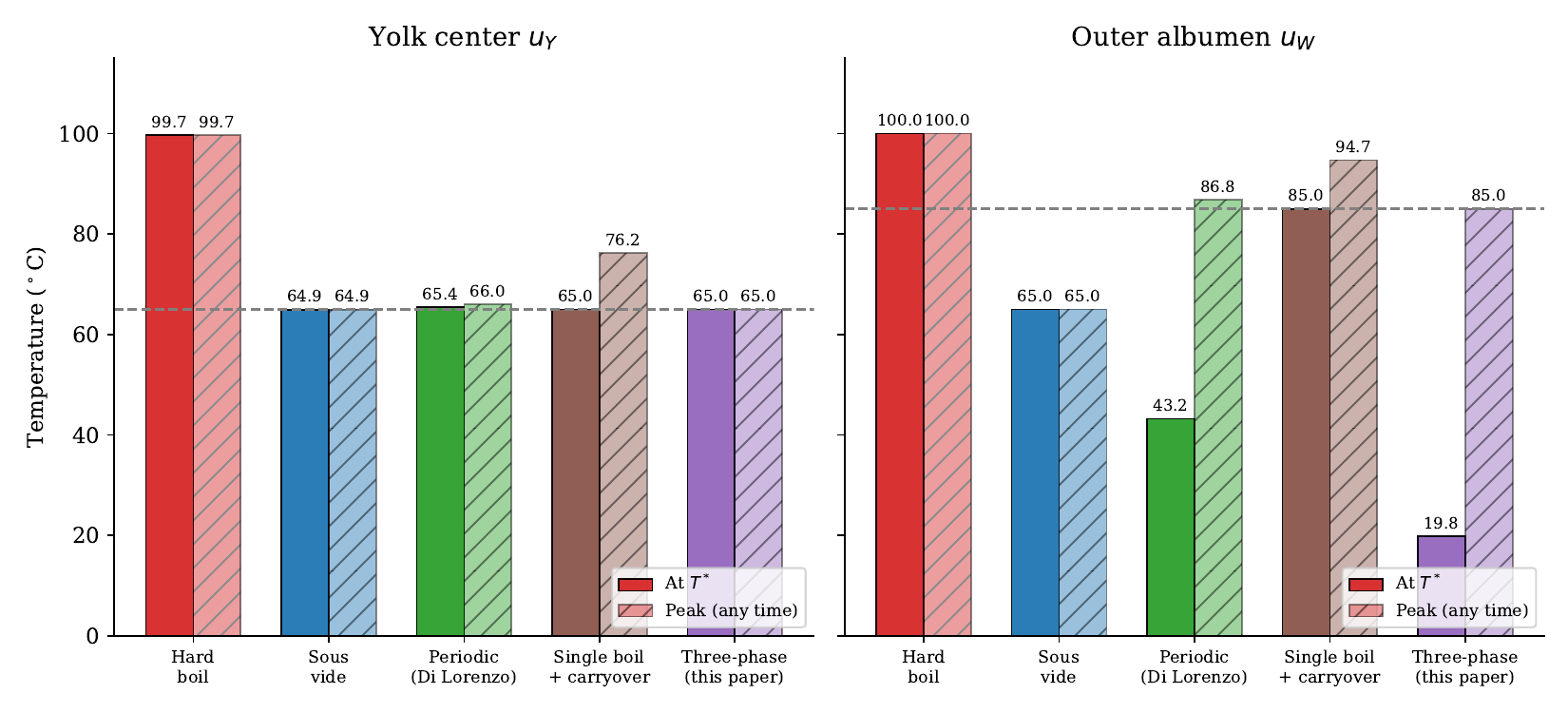}
  \caption{Terminal and peak temperatures for all five protocols.
    Left: yolk center; right: outer albumen.
    Dashed line: target temperature.
    Solid bars show the temperature at $T^*$; hatched bars show the
    maximum temperature reached at any time.
    Only the three-phase protocol (rightmost) reaches both targets
    while keeping every peak at or below its target.}
  \label{fig:comparison}
\end{figure}

Figure~\ref{fig:comparison} summarizes the terminal and peak
temperatures for all five protocols of Table~\ref{tab:comparison}.
The three-phase protocol takes longer than the naive
single-boil-plus-carryover protocol ($20.67$ minutes versus $7.98$
minutes), but this additional time is the price of genuinely never
overcooking either component, rather than merely appearing correct at
one instant; it remains comparable to, and somewhat faster than, the
$32$-minute periodic protocol of \citet{DiLorenzo2025}, while also
being more reliable, since it requires only two changes in water
temperature (sous-vide to boiling, boiling to ice bath) rather than
sixteen alternations, and tolerates small timing errors far better:
because $u_Y$ approaches its peak gradually near $T^*$, the rate of
change of $u_Y$ at $T^*$ is much smaller than in the single-boil
protocol, so a brief delay in removing the egg causes a much smaller
temperature error.

\section{Conclusion}
\label{sec:conclusion}

We solved the two-domain spherical heat equation governing egg
cooking by taking its Laplace transform in time, which reduces the
spatial problem in each domain to an elementary ordinary differential
equation and the boundary/interface conditions to a $3\times3$
algebraic system.
Inverting the resulting transform numerically via Talbot's method,
validated throughout against an independently coded finite-difference
solver, gives a fast and reliable way to evaluate the egg's
temperature under any piecewise-constant cooking protocol.

Using this tool, we found that the periodic protocol of
\citet{DiLorenzo2025} --- alternating boiling and lukewarm water every
two minutes for 32 minutes --- misses both temperature targets at its
terminal time, ending with the yolk hotter than the albumen, and
moreover repeatedly drives the albumen well past its target during the
protocol's high-temperature phases.
A naive single boiling phase, even when corrected by a brief carryover
rest, suffers from the same defect: matching both temperatures exactly
at one instant does not prevent either temperature from having already
overshot its target earlier, nor does it provide any safety margin
against small timing errors.
We showed that avoiding overshoot entirely requires three phases at
different bath temperatures: a sous-vide pre-soak at exactly the
yolk's target temperature, which cannot overshoot by construction; a
short boil to bring the albumen toward its target; and an ice-water
bath that arrests the albumen's overshoot while the heat already
absorbed by the egg continues to bring the yolk to its target.  The
resulting protocol reaches both targets in $20.67$ minutes while never
exceeding either, comparable in total time to the periodic protocol of
\citet{DiLorenzo2025} but considerably more reliable, requiring only
two changes in water temperature rather than sixteen and tolerating
small timing errors far better, since the yolk's temperature changes
only slowly in the neighborhood of $T^*$.

Several extensions merit future study.
First, the thermal diffusivities of yolk and albumen change as
proteins denature, suggesting a nonlinear model in which $\aY$ and
$\aW$ depend on temperature; the present linear analysis would serve
as the leading-order approximation.
Second, the oblate spheroidal geometry of real eggs introduces
anisotropy not captured by the spherical model.
Third, the same Laplace-transform-plus-numerical-inversion strategy
extends naturally to other layered thermal processing problems ---
the tempering of chocolate, the heat treatment of composite spheres,
and the thermal sterilization of heterogeneous food products ---
wherever separate temperature targets must be met in concentric
domains without overshoot, and the underlying eigenfunction expansion
converges too slowly at the times of interest to be evaluated
directly.

We conclude that the mathematical analysis of the soft-boiled egg, far
from being a mere pedagogical curiosity, yields genuinely actionable
insight: a three-phase protocol that, unlike the alternatives
considered here, reaches both temperature targets while never
overcooking either component along the way.

\bibliography{egg-references}

\end{document}